\documentclass[a4paper,12pt]{article}
\usepackage[T1]{fontenc}
\usepackage{latexsym}
\usepackage{amssymb}
\usepackage{graphics}
\newtheorem{thm}{Theorem}
\newtheorem{lem}[thm]{Lemma}
\newtheorem{prop}[thm]{Proposition}
\newtheorem{cor}[thm]{Corollary}

\newcommand{\R}{\mathbb{R}}
\newcommand{\Z}{\mathbb{Z}}
\newcommand{\N}{\mathbb{N}}
\newcommand{\E}{\mathbb{E}}

\title{\textbf{Maximal injective and mixing masas in group factors}}
\author{Paul Jolissaint}

\begin{document}

\maketitle

\begin{abstract}
We present families of pairs of finite von Neumann algebras $A\subset M$ where $A$ is a maximal injective masa in the type $\mathrm{II}_1$ factor $M$ 
with separable predual. Our results 
make use of the strong mixing and the asymptotic orthogonality properties of $A$ in $M$.  
Our method is directly borrowed from S. Popa's original result where he proved that if $G$ is a non abelian free group and if $a$ is one of its 
generators, then the von Neumann algebra generated by $a$ is maximal injective in the factor $L(G)$. Our results apply to pairs $H<G$ where $H$ is an 
infinite abelian subgroup of a
suitable amalgamated product group $G$. 
\par\vspace{3mm}\noindent
\end{abstract}

\section{Introduction}

Inspired by \cite{PopaOA} and \cite{RSS}, we introduced in \cite{JS} the notions of weakly mixing and strongly mixing masas in finite 
von Neumann algebras, and we presented several families of examples coming essentially from pairs of groups. The purpose of the present paper 
is to use these mixing properties in order to give examples of masas that are maximal injective in the ambiant factor. In fact, we will present a 
slight generalization of the main result in S. Popa's article \cite{Pop}, and our point of view is very similar to the exposition of Popa's theorem as 
presented in A. Sinclair's and R. Smith's monograph \cite{SS}. 

Before stating our main results, let us recall our notations. In this article, $M$ will always denote a finite von Neumann algebra with separable 
predual (equivalently, it admits a countable dense set for the strong operator topology), and $\tau$ will be some fixed finite, faithful, normal, 
normalized trace on $M$. If $B$ is a unital von Neumann subalgebra of $M$, then $\E_B$ denotes the $\tau$-preserving conditional expectation onto $B$, 
and we let $M\ominus B$ be the set of all $x\in M$ such that $\E_B(x)=0$; equivalently, it is the set of all $x$ that satisfy $\tau(xb)=0$ for all 
$b\in B$. For future use, we observe that, for all $x,y\in M$ and every $u\in B$, one has:
$$
(*)\quad \E_B((x-\E_B(x))u(y-\E_B(y)))=\E_B(xuy)-\E_B(x)u\E_B(y).
$$ 
If $\Phi$ is a linear map on $M$, we set 
$$
\Vert\Phi\Vert_{\infty,2}=\sup_{x\in M,\Vert x\Vert\leq 1}\Vert \Phi(x)\Vert_2.
$$
If $\omega$ is a free ultrafilter on $\N$, then $M^\omega$ denotes the associated ultrapower algebra; $M$ embeds into $M^\omega$ in a natural way, and, 
if $M$ is a factor, we say that it has \textit{Property} $\Gamma$ of Murray and von Neumann if the relative commutant $M'\cap M^\omega$ is non-trivial. 
If it is the case, it is automatically diffuse. If $M$ does not have Property $\Gamma$, it is called a \textit{full} factor.
For all of this, see for instance Appendix A in \cite{SS}.

Let $G$ be a countable group. We denote by $L(G)$ the von Neumann algebra generated by the left regular representation of $G$ on $\ell^2(G)$; we denote 
simply by $g\xi$ the action of $g\in G$ on $\xi\in\ell^2(G)$ defined by $(g\xi)(h)=\xi(g^{-1}h)$ for every $h\in G$.
As is well known, $L(G)$ is a finite von Neumann algebra, every $x\in L(G)$ has a unique Fourier expansion $\sum_{g\in G}x(g)g$ which converges in the 
$\Vert\cdot\Vert_2$-sense and $\sum_g|x(g)|^2=\Vert x\Vert_2^2$.
Furthermore, $L(G)$ is a factor if and only if $G$ is an ICC group. If $H$ is a subgroup of $G$, then its associated von Neumann algebra $L(H)$ embeds 
into $L(G)$ by setting $x(g)=0$ for all $g\in G\setminus H$ if $x\in L(H)$.

Let now $M$ be a type $\textrm{II}_1$ factor with separable predual and let $A$ be a unital, abelian von Neumann subalgebra of $M$.
Following \cite{JS} and \cite{CFM}, 
we say that $A$ is \textbf{weakly mixing} in $M$ if there exists a sequence of unitary operators $(u_n)\subset A$ such that 
$$
(**)\quad \lim_{n\to\infty}\Vert \E_A(xu_ny)\Vert_2=0\quad \forall x,y\in M\ominus A.
$$
Note that by $(*)$, the latter is equivalent to $\lim_n\Vert \E_A(xu_ny)-\E_A(x)u_n\E_A(y)\Vert_2=0$ for all $x,y\in M$.
We say that $A$ is \textbf{strongly mixing} in $M$ if
$(**)$ holds 
for all sequences of unitary operators $(u_n)\subset A$ such that $\lim_{n\to\infty}u_n=0$ in the weak operator topology. If $G$ is a countable ICC 
group and if $H$ is an abelian subgroup of $G$, if we set $M=L(G)$ and $A=L(H)$,
then it follows from \cite{JS} that $A$ is weakly mixing in $M$ if and only if the pair $H<G$ satisfies the so-called \textit{condition (SS)}: 
for every finite subset $C\subset G\setminus H$, there exists $h\in H$ such that $g_1hg_2\notin H$ for all $g_1,g_2\in C$. Similarly, $A$ is strongly 
mixing in $M$ if and only if the pair $H<G$ satisfies the so-called \textit{condition (ST)}: for every finite subset $C\subset G\setminus H$, there 
exists a finite subset $E\subset H$ such that $g_1hg_2\notin H$ for all $h\in H\setminus E$ and all $g_1,g_2\in C$. 
\par\vspace{3mm}\noindent
\textbf{Remarks.} (1) Let $A$ be a masa in a type $\textrm{II}_1$ factor with separable predual $M$. The main theorem of \cite{SSWW} states that $A$ is 
a singular masa if and only if it is weakly mixing in $M$. See also Theorem 11.1.2 of \cite{SS}.\\
(2) Let $H<G$ be a pair of groups as above. It is easy to see that it satisfies condition (ST) if and only if, for every
$g\in G\setminus H$, the intersection $H\cap gHg^{-1}$ is a finite group. In particular, when the intersection is trivial, $H$ is said to be 
\textit{malnormal} in $G$. Thus, if a subgroup $H$ of a group $G$ which satisfies condition (ST) is also called \emph{almost malnormal} in $G$.

\par\vspace{3mm}
Inspired by \cite{Pop}, the authors of \cite{CFRW} introduced the following property for a pair $A\subset M$ where
$A$ is abelian and $M$ is a type $\mathrm{II}_1$ factor: one
says that $A$ has the \textit{asymptotic orthogonality property} if there is a free ultrafilter $\omega$ on $\N$ such that $x^{(1)}y_1\perp y_2x^{(2)}$ 
in $L^2(M^\omega)$ whenever $x^{(1)},x^{(2)}\in A'\cap M^\omega$ with $\E_{A^\omega}(x^{(i)}))=0$, and $y_1,y_2\in M$ with $\E_A(y_i)=0$ for $i=1,2$. 
Then the authors of \cite{CFRW} prove in Corollary 2.3 that if $A$ is a singular masa which has the asymptotic orthogonality property, then it is 
maximal injective in $M$. 

As we will see, 
strongly mixing masas provide a central decomposition of intermediate algebras that strengthens maximal injectivity. Our first result extends to 
arbitrary pairs $A\subset M$ Theorem 14.2.1 of \cite{SS} which was stated for group algebras; our proof differs partly from that in \cite{SS}. 

\begin{thm}
Let $M$ be a type $\mathrm{II}_1$ factor with separable predual and let $1\in A\subset M$ be a strongly mixing abelian von Neumann subalgebra of $M$. 
If $N$ is a von Neumann subalgebra of $M$ which contains $A$, then there exists a partition of the unity $(e_k)_{k\geq 0}$ in the center $Z$ of $N$ such 
that $Ne_0=Ae_0$ and, for every $k\geq 1$ such that $e_k\not=0$, $Ne_k$ is a type $\mathrm{II}_1$ factor and $(N'\cap A^\omega)e_k$ has a non-zero atomic part.
\end{thm}
\par\vspace{3mm}
The strong mixing assumption in Theorem 1 is essential, as was kindly communicated to us by S. White. Indeed, let $A_0\subset M_0$ be a weakly, but not strongly mixing masa in a type $\mathrm{II}_1$ factor. Then $A=A_0\overline{\otimes}A_0$ is obviously a weakly mixing masa in the factor $M=M_0\overline{\otimes}M_0$, but taking $N=A_0\overline{\otimes}M_0$, we have $Z=A_0\otimes 1$ which has no atoms.
\par\vspace{3mm}
The second result is essentially Theorem 14.2.5 of \cite{SS} where it was stated in the special case of the free group factors.
We recall it for the sake of completeness and future use.
\begin{thm}
Let $A$ be a strongly mixing masa that satisfies the asymptotic orthogonality property
in a type $\mathrm{II}_1$ factor $M$ with separable predual, let $N$ be an intermediate von Neumann subalgebra and let $(e_k)_{k\geq 0}\subset Z(N)$ be the corresponding partition of the unity as in Theorem 1. Then 
for every $k\geq 1$ such that $e_k\not=0$, $Ne_k$ is a full type $\mathrm{II}_1$ factor. In particular, $A$ is a maximal injective subalgebra of $M$.
\end{thm}

As will be seen, pairs of groups can provide examples of such algebras.
Thus, for the rest of the present section, let $G$ be an ICC countable group and let $H$ be an abelian subgroup of $G$. Put $M=L(G)$ and $A=L(H)$. 
We assume that the pair $H<G$ satisfies the following hypotheses:\\
There exists a sequence $(W_m)_{m\geq 1}$ of subsets of $G\setminus H$ such that
\begin{enumerate}
	\item [(H1)] $W_m\subset W_{m+1}$ for every $m\geq 1$ and $\bigcup_m W_m=G\setminus H$ ;
	\item [(H2)] if $V_m$ denotes the complementary set of $W_m\cup W_m^{-1}$ in $G\setminus H$, then for all $g_1,g_2\in G\setminus H$, there exists a
	positive integer $m_1=m_1(g_1,g_2)$ such that $g_1V_m\cap V_mg_2=\emptyset$ for every $m>m_1$;
	\item [(H3)] there exist an integer $m_0>0$ and an element $h\in H$ such that, for every $m>m_0$, one can find an integer $i_m>0$ such that $h^iW_mh^{-i}\cap W_m=\emptyset$ for every $i\geq i_m$.
\end{enumerate}

\begin{thm}
Let $H<G$ be a pair of groups such that $G$ is countable and ICC, $H$ is abelian, and assume that the pair satisfies condition (ST) on the one hand, and conditions (H1), (H2) and (H3) on the other hand. Then $L(H)$ is a strongly mixing masa that satisfies the asymptotic orthogonality property in $L(G)$. Thus, 
if $N$ is a von Neumann subalgebra of $M$ which contains $A$, then there exists a partition of the unity $(e_k)_{k\geq 0}$ in the center of $N$ such that $Ne_0=Ae_0$ and, for every $k\geq 1$ such that $e_k\not=0$, $Ne_k$ is a full type $\mathrm{II}_1$ factor. In particular, $A$ is a maximal injective subalgebra of $M$.
\end{thm}

\par\vspace{3mm}
The next section is devoted to the proof of Theorems 1, 2 and 3. In Section 3, we provide a family of examples of pairs $H<G$ that satisfy all conditions of Theorem 2, hence which gives examples of maximal injective masas in group factors; it is given by amalgamated products $G=H*_ZK$ where $H$ is infinite and abelian, $Z$ is finite and different from $K$, and $G$ is ICC. 

\par\vspace{3mm}
\emph{Acknowledgements.} I am grateful to Stuart White 
for helpful
comments and for having detected a mistake in the first version of Theorem 1.

\section{Proofs of the main results}

Let $M$ be a type $\mathrm{II}_1$ factor with separable predual and let $A$ be a masa in $M$. Before proving Theorem 1, we present an auxiliary result of independent interest. See also Section 4 in \cite{CFM} for related results.

\begin{prop}
Let $N$ be a finite von Neumanna algebra with separable predual and let $A\subset N$ be a strongly mixing abelian von Neumann algebra in $N$. Then:
\begin{enumerate}
	\item [(1)] the von Neumann algebra $N'\cap A^\omega$ has a non-zero atomic part;
	\item [(2)] for every diffuse von Neumann subalgebra $1\in B\subset A$ and for every unitary $u\in U(N)$, one has 
	$$
\Vert \E_B-\E_{uBu^*}\Vert_{\infty,2}\geq \Vert u-\E_{A}(u)\Vert_2.
$$
In particular, $\mathcal{N}_N(B)''\subset A$.
	
\end{enumerate}

\end{prop}
\textit{Proof.} (1) We assume that $N'\cap A^\omega$ is diffuse
and we will get a contradiction. Choose a Haar unitary $u\in U(N'\cap A^\omega)$; this means that $\tau_\omega(u^k)=0$ for all integers $k\not=0$. The algebra $N$ being separable with respect to the $\Vert\cdot\Vert_2$-topology, choose an increasing sequence of finite subsets 
$$
E_1\subset E_2\subset \ldots \{x\in N\ominus A:\Vert x\Vert\leq 1\}
$$
so that $\bigcup_n E_n$ is dense in $\{x\in N\ominus A:\Vert x\Vert\leq 1\}$. Similarly, choose an increasing sequence of finite subsets $F_1\subset F_2\ldots$ of the unit ball of $A$ which is $\Vert\cdot\Vert_2$-dense. Let us write $u=[(u_n)_n]$ with $u_n\in U(A)$. The unitaries $u^k$ being pairwise orthogonal in $A^\omega$, one has, by Parceval's inequality
$$
\sum_{k\in\Z}|\tau_\omega(xu^k)|^2\leq\Vert x\Vert_2^2
$$
hence
$$
\lim_{|k|\to\infty}\tau_\omega(xu^k)=0
$$
for every $x\in A^\omega$. Thus, for every integer $n\geq 1$, there exists an integer $k_n>0$ such that
$$
\max_{a\in F_n}|\tau_\omega(au^k)|<\frac{1}{2n}\quad \forall |k|\geq k_n.
$$
As moreover $xu^k=u^kx$ for every $x\in N$ and every integer $k$, we infer that, for every $n>0$, there exists a positive integer $\ell_n$ such that
$$
|\tau(au_{\ell_n}^{k_n})|<\frac{1}{n}\quad \forall a\in F_n\quad\textrm{and}\quad 
\Vert u_{\ell_n}^{k_n}xu_{\ell_n}^{-k_n}-x\Vert_2<\frac{1}{n}\quad\forall x\in E_n.
$$
Thus the sequence $(u_{\ell_n}^{k_n})\subset U(A)$ converges to 0 in the weak operator topology, and we get for every $x\in N\ominus A$ 
$$
\lim_{n\to\infty}\Vert \E_A(x^*u_{\ell_n}^{k_n}x)\Vert_2=0.
$$
However, if $x\not=0$ is orthogonal to $A$, we have
$$
0<\Vert \E_A(x^*x)\Vert_2\leq\Vert \E_A(x^*(x-u^{k_n}_{\ell_n}xu^{-k_n}_{\ell_n}))\Vert_2+\Vert 
\E_A(x^*u^{k_n}_{\ell_n}x)\Vert_2\to 0,
$$
which gives the desired contradiction.\\
(2) Fix a diffuse von Neumann subalgebra $B$ of $A$, 
a unitary operator $u\in U(N)$, and let us consider $x=u^*-\E_{A}(u^*)$ and $y=u-\E_{A}(u)\in N\ominus A$ and let $\varepsilon>0$. One has for every $v\in U(B)$:
$$
\E_{B}(xvy)  = 
\E_{B}(u^*vu)-\E_{B}(\E_{A}(u^*)v\E_{A}(u))
$$
As $B$ is diffuse, there exists a sequence of unitaries $(v_n)\subset U(B)$ which converges to 0 with respect to the weak operator topology.
Since $A$ is strongly mixing in $N$, there exists a positive integer $n$ such that $\Vert \E_{A}(xv_ny)\Vert_2\leq \varepsilon$. As $\E_B=\E_B\E_{A}$, we have 
$\Vert \E_{B}(xv_ny)\Vert_2\leq \varepsilon$ as well.
The above computation gives
$$
\Vert \E_{B}(u^*v_nu)\Vert_2\leq \Vert \E_{B}(\E_{A}(u^*)v_n\E_{A}(u))\Vert +\varepsilon
\leq \Vert \E_{A}(u)\Vert_2 +\varepsilon.
$$
We get then:
\begin{eqnarray*}
\Vert \E_{B}-\E_{uBu^*}\Vert_{\infty,2}^2 & \geq &
\Vert v_n-\E_{uBu^*}(v_n)\Vert_2^2\\
&=& \Vert u^*v_nu-\E_{B}(u^*v_nu)\Vert_2^2\\
&=&
1-\Vert \E_{B}(u^*v_nu)\Vert_2^2\\
&\geq &
1-(\Vert \E_{A}(u)\Vert_2+\varepsilon)^2\\
&=&
1-\Vert \E_{A}(u)\Vert_2^2-2\varepsilon\Vert \E_{A}(u)\Vert_2-\varepsilon^2\\
&=&
\Vert u-\E_{A}(u)\Vert_2^2-2\varepsilon\Vert \E_{A}(u)\Vert_2-\varepsilon^2.
\end{eqnarray*} 
As $\varepsilon$ is arbitrary, we get the conclusion.
\hfill $\square$
\par\vspace{3mm}\noindent
\textbf{Remark.} Let $A$ be a masa in the type $\mathrm{II}_1$ factor $M$. Let us say that it satisfies condition (D) if, for every diffuse von Neumann subalgebra $1\in B\subset A$, the normalizer $\mathcal{N}_M(B)$ is equal to the unitary group $U(A)$. Statement (2) of Proposition 4 above shows that if $A$ is strongly mixing in $M$ then it satisfies condition (D), and one may ask whether the converse holds true. In fact, here is a partial positive answer: suppose that $G$ is an ICC group and that $H$ is an infinite abelian subgroup of $G$. Set $A=L(H)$ and $M=L(G)$. Then if $A\subset M$ satisfies condition (D), $A$ is strongly mixing in $M$. Indeed, suppose that it is not the case; by Theorem 3.5 in \cite{JS}, the pair $H<G$ would not satisfy condition (ST), hence
there exists an element $g\in G\setminus H$ such that $H\cap gHg^{-1}$ is an infinite subgroup of $H$. 
Set $B=L(H\cap gHg^{-1})$, which is a diffuse subalgebra of $A$. Then $g\in\mathcal{N}_M(B)$ but it does not belong to $U(A)$.

\par\vspace{3mm}\noindent
\textit{Proof of Theorem 1.} Our proof is strongly inspired by that of theorem 14.2.1 in \cite{SS}.
Let $e_0\in Z$ be the largest projection such that $Ne_0=Ae_0$ and set $e=1-e_0$. We will prove that $eZ$ is atomic. Assume on the contrary that it is not. Then there exists a projection $q\in Z$ such that $0\not=q\leq e$ and $qZ$ is a diffuse algebra and $(e-q)Z$ is atomic. Set $B=(1-q)A+qZ$. As $A$ is maximal abelian and since $Z$ commutes with $A$, we have $Z\subset A$ hence $B\subset A$ with $qB=qZ$ diffuse, hence $B$ is diffuse as well. By the previous proposition, one has $\mathcal{N}_M(B)''\subset A$, and in particular $\quad q(B'\cap M)\subset A$.
As in the proof of Theorem 14.2.1 in \cite{SS}, we get
$$
qN\subset q(Z'\cap M)=q(B'\cap M)\subset qA,
$$
hence $qA=qN$ which is diffuse, and this contradicts maximality of $e_0$ since $q\leq 1-e_0$. 
\par

Finally, as observed in \cite{CFM}, for every non-zero projection $e\in A$, the masa $Ae$ is still strongly mixing in $eMe$. Hence $Ae_k$ is strongly mixing in $e_kMe_k$, and we deduce that $(N'\cap A^ \omega)e_k$ has atoms by the previous proposition.
\hfill $\square$

\par\vspace{3mm}\noindent
\textit{Proof of Theorem 2.} We reproduce the proof of Theorem 14.2.5 in \cite{SS} for the reader's convenience.
If $A\subset N\subset M$ are as in Theorem 2, let $(e_k)_{k\geq 0}\subset Z(N)$ be the partition of the unity provided by Theorem 1: $Ne_0=Ae_0$, and for every $k>0$ such that $e_k\not=0$, the von Neumann algebra
$Ne_k$ is a type $\mathrm{II}_1$ factor such that the von Neumann algebra $(N'\cap A^\omega)e_k$ has a non-zero atomic part. If there is some $k>0$ such that $Ne_k$ has Property $\Gamma$, then the relative commutant $(Ne_k)'\cap (Ne_k)^\omega$ is diffuse, hence $(N'\cap A^\omega)e_k\subsetneq (Ne_k)'\cap (Ne_k)^\omega$. As in the proof of Theorem 14.2.5 of \cite{SS}, we choose some non-zero $x\in (Ne_k)'\cap (Ne_k)^\omega$ such that $\E_{A^ \omega}(x)=0$ and some unitary $w\in Ne_k$ such that $\E_A(w)=0$. By the asymptotic orthogonality property of $A\subset M$ applied to $x^{(1)}=x^{(2)}=x$ and $y_1=y_2=w$,
we get
$$
2\Vert x\Vert_2^2=\Vert wx\Vert_2^2+\Vert xw\Vert_2^2=\Vert wx-xw\Vert_2^2=0
$$
which is a contradiction.
\hfill $\square$

\par\vspace{3mm}
From now on, we consider a countable, ICC group $G$ and an abelian, infinite subgroup $H$ of $G$ and we assume that the pair $H<G$ satisfies the three conditions (H1) to (H3) in Section 1, and we set as before $A=L(H)\subset M=L(G)$. For convenience, we recall some notations from \cite{SS}: 
If $W\subset G$, let $p_W$ be the orthogonal projection of $\ell^2(G)$ onto the subspace $\ell^2(W)$:
$$
p_W(x)=\sum_{g\in W}x(g)g\quad\textrm{if}\quad x=\sum_{g\in G}x(g)g.
$$
We remind the reader that for all $V,W\subset G$, one has
$p_Vp_W=p_{V\cap W}$ (thus in particular, $p_Vp_W=0$ if $V$ and $W$ are disjoint),
that
$$
p_{gWg^{-1}}(x)=gp_W(g^{-1}xg)g^{-1}\quad\textrm{and}\quad p_W(x)^*=p_{W^{-1}}(x^*)
$$
for all $g\in G$ and $x\in\ell^2(G)$, and that, if $V\subset W$, then $\Vert p_V(x)\Vert_2\leq\Vert p_W(x)\Vert_2$ for every $x$.

\par\vspace{3mm}
The proof of Theorem 3 follows immediately from part (iii) of the
following lemma whose proof is similar to those of Lemmas 14.2.3 and 14.2.4 in \cite{SS}. However, we give a proof for the sake of completeness.

\begin{lem}
Let $G$, $H$ satisfy conditions (H1) and (H2), and let $(W_m)_{m\geq 1}$ be the corresponding sequence of subsets of $G\setminus H$. Assume also that it satisfies the following weaker variant of condition (H3):\\
(H3') there exists an integer $m_0>0$ such that, for every $m>m_0$, one can find elements $h_{1,m},\ldots,h_{n_m,m}\in H$ such that 
	$n_m\to\infty$ as $m\to\infty$ and
	$$
	h_{i,m}W_mh_{i,m}^{-1}\cap h_{j,m}W_mh_{j,m}^{-1}=\emptyset\quad\forall i\not=j.
	$$
\begin{enumerate}
\item [(i)] Let $\varepsilon>0$, let $m>m_0$ and $h_{1,m},\ldots,h_{n_m,m}$ 
be as in condition (H3') and let 
$x\in\ell^2(G)$, $\Vert x\Vert_2\leq 1$, be such that
$$
\Vert h_{j,m}xh_{j,m}^{-1}-x\Vert_2\leq \varepsilon\quad \forall j=1,\ldots,n_m.
$$
Then 
$$
\Vert p_{W_m\cup W_m^{-1}}(x)\Vert_2^2\leq 4(\varepsilon^2+n_m^{-1}).
$$
\item [(ii)] If $y\in\ell^2(G)$ is such that $\E_A(y)=0$, then
	$$
	\lim_{m\to\infty}\Vert y-p_{W_m}(y)\Vert_2=0.
	$$
\item [(iii)] The abelian algebra $A$ satisfies the asymptotic orthogonality property in $M$, namely,
let $\omega$ be a free ultrafilter on $\N$, $x^{(1)},x^{(2)}\in A'\cap M^\omega$ and $y_1,y_2\in M$ be such that $\E_{A^\omega}(x^{(j)})=\E_A(y_j)=0$ for $j=1,2$.
Then $y_1x^{(1)}\perp x^{(2)}y_2$ in $M^\omega$ and
$$
\Vert y_1x^{(1)}-x^{(2)}y_2\Vert_{\omega,2}^2=\Vert y_1x^{(1)}\Vert_{\omega,2}^2+\Vert x^{(2)}y_2\Vert_{\omega,2}^2.
$$	
\end{enumerate}
\end{lem}
\textit{Proof.} (i) Using $(\alpha+\beta)^2\leq 2(\alpha^2+\beta^2)$ for arbitrary real numbers $\alpha$ and $\beta$, we get for every $m>m_0$ and every $1\leq j\leq n_m$:
\begin{eqnarray*}
\Vert p_{W_m}(x)\Vert_2^2 & \leq &
(\Vert p_{W_m}(x-h_{j,m}^{-1}xh_{j,m})\Vert_2+\Vert p_{W_m}(h_{j,m}^{-1}xh_{j,m})\Vert_2)^2\\
&\leq &
2\Vert x-h_{j,m}^{-1}xh_{j,m}\Vert_2^2+2\Vert p_{W_m}(h_{j,m}^{-1}xh_{j,m})\Vert_2^2\\
&\leq &
2\varepsilon^2+2\Vert p_{h_{j,m}W_mh_{j,m}^{-1}}(x)\Vert_2^2.
\end{eqnarray*} 
Using condition (H3'), $p_{h_{j,m}W_mh_{j,m}^{-1}}(x)$ is orthogonal to $p_{h_{i,m}W_mh_{i,m}^{-1}}(x)$ for all $i\not=j$. Summing over 
$1\leq j\leq n_m$, we get:
\begin{eqnarray*}
n_m\Vert p_{W_m}(x)\Vert_2^2 &\leq &
2n_m\varepsilon^2+2\sum_{j=1}^{n_m}\Vert p_{h_{j,m}W_mh_{j,m}^{-1}}(x)\Vert_2^2\\
&\leq &
2n_m\varepsilon^2+2\Vert\sum_{j=1}^{n_m} p_{h_{j,m}W_mh_{j,m}^{-1}}(x)\Vert_2^2\\
&\leq &
2n_m\varepsilon^2+2.
\end{eqnarray*}
Hence
$$
\Vert p_{W_m}(x)\Vert_2^2 \leq 2(\varepsilon^2+n_m^{-1}).
$$
As $x^*$ satisfies the same conditions as $x$, using  $p_{W_m}(x^*)=p_{W_m^{-1}}(x)^*$, we have
\begin{eqnarray*}
\Vert p_{W_m\cup W_m^{-1}}(x)\Vert_2^2 
& = &
\Vert p_{W_m}(x)\Vert_2^2+\Vert p_{W_m^{-1}\setminus W_m}(x)\Vert_2^2\\
&\leq &
\Vert p_{W_m}(x)\Vert_2^2+\Vert p_{W_m^{-1}}(x)\Vert_2^2\\
&\leq &
4(\varepsilon^2+n_m^{-1}).
\end{eqnarray*}
This proves claim (i).\\
Claim (ii) follows immediately from condition (H1).\\
Let us prove claim (iii): We assume that $\Vert x^{(i)}\Vert,\Vert y_i\Vert\leq 1$ for $i=1,2$. Furthermore, $x^{(i)}=[(x^{(i)}_r)_r]$, and, replacing $x^{(i)}_r$ by $x^{(i)}_r-\E_A(x^{(i)}_r)$, 
we assume that $\E_A(x^{(i)}_r)=0$ and $\Vert x^{(i)}_r\Vert\leq 1$ for every $r\in \N$ and $i=1,2$. Suppose first that $y_j=g_j\in G\setminus H$ for $j=1,2$.
\par 
Let $\varepsilon>0$ be fixed; we will show that
$$
|\tau_\omega(g_1x^{(1)}[x^{(2)}g_2]^*)|\leq 6\varepsilon.
$$
Let us choose $m>0$ large enough so that $g_1V_m\cap V_mg_2=\emptyset$ and that
$n_m^{-1}\leq\varepsilon^2$. Thus, $g_1p_{V_m}(x^{(1)}_r)\perp p_{V_m}(x^{(2)}_r)g_2$ for every $r$.
\par 
Next set
$$
T=\{r\in\N:\Vert h_{j,m}x^{(i)}_r-x^{(i)}_rh_{j,m}\Vert_2\leq\varepsilon,\ \forall 1\leq j\leq n_m\ \mathrm{and}\ i=1,2\}
$$
which belongs to $\omega$ since each $x^{(i)}\in A'\cap M^\omega$. By part (ii), we have for $r\in T$ and $i=1,2$,
$$
\Vert x^{(i)}_r-p_{V_{m}}(x^{(i)}_r)\Vert_2^2=\Vert p_{W_{m}\cup W_{m}^{-1}}(x^{(i)}_r)\Vert_2^2\leq 4\left(\varepsilon^2+\frac{1}{n_m}\right)\leq 8\varepsilon^2.
$$
(We used $p_H(x^{(i)}_r)=0$, hence $x^{(i)}_r=p_{V_{m}}(x^{(i)}_r)+p_{W_{m}\cup W_{m}^{-1}}(x^{(i)}_r)$.)
\par
For the same values of $r$, we get:
\begin{eqnarray*}
|\tau(g_1x^{(1)}_r[x^{(2)}_rg_2]^*)| 
&\leq &
|\tau(g_1(x^{(1)}_r-p_{V_m}(x^{(1)}_r))g_2^{-1}x^{(2)*}_r)|
+
|\tau(g_1p_{V_{m}}(x^{(1)}_r)g_2^{-1}x^{(2)*}_r)|\\
&\leq &
\Vert x^{(1)}_r-p_{V_{m}}(x^{(1)}_r)\Vert_2
+|\tau(g_1p_{V_{m}}(x^{(1)}_r)g_2^{-1}(x_r^{(2)*}-p_{V_{m}}(x_r^{(2)*})))|\\
& &
+|\tau(g_1p_{V_{m}}(x^{(1)}_r)g_2^{-1}p_{V_m}(x_r^{(2)*}))|\\
&\leq &
\Vert x^{(1)}_r-p_{V_{m}}(x^{(1)}_r)\Vert_2+\Vert x^{(2)}_r-p_{V_{m}}(x^{(2)}_r)\Vert_2\\
& &
+|\tau([g_1p_{V_{m}}(x^{(1)}_r)][p_{V_{m}}(x^{(2)}_r)g_2]^*)|\\
&\leq & 2\sqrt{8}\varepsilon<6\varepsilon
\end{eqnarray*}
since $\tau([g_1p_{V_{m}}(x^{(1)}_r)][p_{V_{m}}(x^{(2)}_r)g_2]^*)=0$. Thus
$$
T\subset T':=\{r\in\N:|\tau(g_1x^{(1)}_r[x^{(2)}_rg_2]^*)| \leq 6\varepsilon\}
$$
and $T'\in\omega$, hence $|\tau_\omega(g_1x^{(1)}_r[x^{(2)}_rg_2]^*)| \leq 6\varepsilon$. 
\par
Next, using linearity, $\tau_\omega(y_1x^{(1)}[x^{(2)}y_2]^*)=0$ for all $y_1,y_2\in M$ with finite support and such that $\E_A(y_1)=\E_A(y_2)=0$, 
and using density and the same kind of arguments as above, we get $\tau_\omega(y_1x^{(1)}[x^{(2)}y_2]^*)=0$ for arbitrary $y_1,y_2\in M\ominus A$.
\par
Finally, the equality 
$$
\Vert y_1x^{(1)}-x^{(2)}y_2\Vert_{\omega,2}^2=\Vert y_1x^{(1)}\Vert_{\omega,2}^2+\Vert x^{(2)}y_2\Vert_{\omega,2}^2
$$
comes from $y_1x^{(1)}\perp x^{(2)}y_2$.
\hfill $\square$

\section{Examples}

Before discussing our first family of examples, we need to recall some facts on length-functions on groups
taken from \cite{JRD}; a \textit{length-function} on a group $\Gamma$ is a map $\ell:\Gamma\rightarrow\R_+$ satisfying:
\begin{enumerate}
\item [(i)] $\ell(gh)\leq \ell(g)+\ell(h)$ for all $g,h\in\Gamma$;
\item [(ii)] $\ell(g^ {-1})=\ell(g)$ for every $g\in\Gamma$;
\item [(iii)] $\ell(e)=0$, where $e$ denotes the identity in $\Gamma$.
\end{enumerate}
Typical and important examples of length-functions are provided by \textit{word} length-functions in finitely generated groups: if $\Gamma$ is a finitely generated group and if $S$ is a finite, symmetric set of generators of $\Gamma$, then the associated word length-function is defined by
$$
\ell_S(g)=\min\{n\in\N:g=s_1\cdots s_n,\ s_i\in S\}.
$$
If $S'$ is another finite generating set then $\ell_S$ and $\ell_{S'}$ are \textit{equivalent} in the sense that there exist positive numbers $a,a'$ such that
$$
\ell_{S'}(g)\leq a\ell_{S}(g)\quad\textrm{and}\quad \ell_S(g)\leq a'\ell_{S'}(g)
$$
for every $g\in \Gamma$. When the generating set $S$ is fixed, one often write $|g|$ instead of $\ell_S(g)$.
\par\vspace{3mm}
Let now $G=H_1*_Z H_2$ be an amalgamated product where $H_1$ and $H_2$ are finitely generated groups,
$H_1$ is infinite and abelian, $Z$ is a common finite subgroup of $H_1$ and $H_2$, $Z\not=H_2$, and we assume that $G$ is an ICC group.
We choose sets of representatives $R_1\ni e$ and $R_2\ni e$ of left $Z$-cosets in $H_1$ and $H_2$ respectively, and, because $Z$ is a finite group, we choose length-functions $\ell_1$ and $\ell_2$ on $H_1$ and $H_2$ respectively with the following properties
(cf \cite{JRD}, Section 2.2):
\begin{enumerate}
 \item [(a)] $\ell_1$ and $\ell_2$ take integer values and are equivalent to the word length-functions on $H_1$ and $H_2$ respectively;
\item [(b)] for all $z,w\in Z$, $j=1,2$ and all $h\in H_j$, one has $\ell_j(zhw)=\ell_j(h)$;
\item [(c)] $\{h\in H_j:\ell_j(h)=0\}=Z$ for $j=1,2$.
\end{enumerate}
We set hereafter $|h|=\ell_j(h)$ for $j=1,2$ and $h\in H_j$ and we observe that, for every $h\in H_j$, if $h=rz$ denotes its decomposition with $r\in R_j$ and $z\in Z$, then $|h|=|r|$. We recall that every 
$g\in G$ has a unique normal form
$$
g=r_1\cdots r_nz
$$
with $n\geq 0$, $z\in Z$, and, if $n>0$, then $r_j\in R_{i_j}\setminus\{e\}$ and $i_j\not=i_{j+1}$ for every $j=1,\ldots,n-1$. For such a $g$, put
$$
|g|=|r_1|+\cdots+|r_n|;
$$
this defines a length-function on $G$ that is equivalent to the word length-function. 
\par\vspace{3mm}
Thus we define for every $m\geq 1$ : 
$$
 W_m=(H_2\setminus Z)\cup\{g=r_1\cdots r_nz: n\geq 2, r_1\in R_1, |r_1|< m\}.
$$

Let us check that $W_m$ and $V_m=(W_m\cup W_m^{-1})^c$ satisfy conditions (H1) to (H3) of the first section. Indeed, (H1) is obviously satisfied. For
(H2), notice first that the normal form of every $g\in V_m$ is of the following type:
$$
(*)\ g=g_1g_2\cdots g_kz
$$ 
where $k\geq 3$, $g_1,g_k\in R_1$, $|g_1|,|g_k|\geq m$ and $z\in Z$.
Thus,
if $\gamma,\gamma'\in (H_1*_Z H_2)\setminus H_1$ are fixed and if $n>0$ is such that $|\gamma|,|\gamma'|<n$ then for every $m>2n$ and all
$g,g'\in V_m$, the element $\gamma g$ cannot start as in $(*)$ and $g'\gamma'$ cannot end as in $(*)$ either.
Hence they cannot be equal.

Finally, let us check that (H3) holds true.
Since $H_1$ is infinite, abelian and finitely generated, it follows from  the structure of such groups that one can choose an element $h\in H_1\setminus Z$ of infinite order. Moreover, if $i\not=j$ are integers, $h^i$ and $h^j$ belong to 
different cosets $\bmod\ Z$ (otherwise $Z$ would contain an element of infinite order). Hence we can assume that $h^j\in R_1$ for every $j\in\Z$. One also has $\lim_{|j|\to\infty}|h^j|=\infty$. Let $j_m>0$ be large enough so that $|h^j|>2m$ for every $j>j_m$; then one has $h^jW_mh^{-j}\cap W_m=\emptyset$ for all such $j$'s. This proves that the sequence $(W_m)$ satisfies condition (H2).

Thus we get, using the fact that the pair $H_1<G$ satisfies also condition (ST) from Proposition 4.1 of \cite{JS}:

\begin{cor}
Let $G=H_1*_Z H_2$ be an amalgamated product as above. If $L(H_1)\subset N\subset L(G)$ is an intermediate von Neumann algebra, there exists partition of the unity $(e)_{k\geq 0}$ in the center of $N$ such that $L(H_1)e_0=Ne_0$, and, for each $k$ such that $e_k\not=0$, $Ne_k$ is a full factor. In particular, $L(H_1)$ is strongly mixing and maximal injective in $L(G)$.
\end{cor}

Similarly, let $G=K*L$ be a free product group such that $|K|\geq 2$ and $L$ contains an element $\beta$ of order at least 3. Let $\alpha$ be some non-trivial element of $K$ and set $H=\langle \alpha\beta\rangle$. Then, by Corollary 4.5 of \cite{JS}, the pair $H<G$ satisfies condition (ST) and it is easy to see that it satisfies also conditions (H1) to (H3) of Section 1. Thus we get:

\begin{cor}
With $H<G=K*L$ as above, and let $L(H)\subset N\subset L(G)$ is an intermediate von Neumann algebra. Then there exists partition of the unity $(e)_{k\geq 0}$ in the center of $N$ such that $L(H)e_0=Ne_0$, and, for each $k$ such that $e_k\not=0$, $Ne_k$ is a full factor. In particular, $L(H)$ is strongly mixing and maximal injective in $L(G)$.
\end{cor}

\par\vspace{3mm}
The following proposition is straightforward.

\begin{prop}
Let $G_1$ be a countable ICC group and let $H<G_1$ be an infinite abelian subgroup. Assume that $G_1\setminus H$ contains a sequence $(W_m)$ of subsets which satisfies conditions (H1), (H2) and (H3) of Section 1. Let $G_2$ be an arbitrary, at most countable, non-trivial group and let $G=G_1*G_2$ be the corresponding free product. For every $m>0$, let $W'_m$ be the set of reduced words $w=g_1\cdots g_n\in G_1*G_2\setminus H$ such that either $g_1\in W_m$ or $w=h^kg_2\cdots g_n$ with $0\leq |k|<m$ and $g_2\in G_2\setminus\{e\}$. Then the sequence $(W'_m)$ satisfies conditions (H1), (H2) and (H3).
\end{prop}

As a consequence of Proposition 3.7 of \cite{JS}, if $L(H)$ is strongly mixing in $L(G)$, then it is also strongly mixing in the free product factor $L(G_1*G_2)=L(G_1)*L(G_2)$, thus we get:

\begin{cor}
If $H<G_1$ and $G_2$ are as in Proposition 8 and if
$L(H)$ is strongly mixing in $L(G_1)$, then $L(H)$ is maximal injective in $L(G_1)$ and in $L(G_1*G_2)$.
\end{cor}

\par\vspace{1cm}
\bibliographystyle{plain}
\bibliography{max}
\par
\vspace{1cm}
\noindent
\begin{flushright}
     \begin{tabular}{l}
       Universit\'e de Neuch\^atel,\\
       Institut de Math\'emathiques,\\       
       Emile-Argand 11\\
       Case postale 158\\
       CH-2009 Neuch\^atel, Switzerland\\
       \small {paul.jolissaint@unine.ch}
     \end{tabular}
\end{flushright}

\end{document}